\documentclass{article}
\usepackage[english]{babel}

\usepackage[letterpaper,top=2cm,bottom=2cm,left=2cm,right=2cm,marginparwidth=1.75cm]{geometry}

\usepackage{float}
\usepackage{calrsfs}
\usepackage{amssymb}
\usepackage{amsmath}
\usepackage{amsthm}
\usepackage{subcaption}
\usepackage{multirow}
\usepackage{graphicx}
\usepackage{amsfonts}
\usepackage{mathtools}
\usepackage{mathrsfs}
\usepackage{placeins}
\usepackage{comment}
\usepackage{bbm}
\usepackage{dsfont}
\usepackage[export]{adjustbox}
\usepackage{overpic}

\usepackage{diagbox}
\usepackage{algorithm}
\usepackage{algpseudocode}

\usepackage{booktabs}

\usepackage{empheq}
\usepackage[mathscr]{euscript}
\usepackage[dvipsnames,x11names]{xcolor}
\usepackage[colorlinks=true, allcolors=OliveGreen]{hyperref}
\usepackage{amsmath}   
\usepackage{tikz}
\usepackage{pgfplots}
\usetikzlibrary{external}
\usepackage{authblk}
\usepackage{array}
\newcolumntype{C}{>{\centering\arraybackslash}p{2.5cm}}
\usepackage{nicematrix}
\usepackage{caption}
\usepackage{subcaption}
\usepackage{tikz}

\usetikzlibrary{patterns,tikzmark}
\usetikzlibrary{matrix,decorations.pathreplacing,calc}
\usepackage{hf-tikz}

\usepackage{orcidlink}
\usepackage{csvsimple}
\usepackage{siunitx}
\usepackage{tikz} 

\usetikzlibrary{arrows.meta, positioning, shapes.geometric}
\usepackage{pgfplots}
\pgfplotsset{compat=newest} 
\usepgfplotslibrary{units} 
\usepackage{tikz}
\usetikzlibrary{matrix,decorations.pathreplacing}
\pgfkeys{/pgfplots/legend entry/.code=\addlegendentry{#1}}

\usepackage{pdfpages}

\newcommand{\averagel}{\{\!\!\{}
\newcommand{\averager}{\}\!\!\}}

\newcommand{\jumpl}{[\![}
\newcommand{\jumpr}{]\!]}

\usepackage[export]{adjustbox}
\usepackage{lipsum}                     
\usepackage{xargs}                      
\usepackage[colorinlistoftodos,prependcaption]{todonotes}
\newcommandx{\unsure}[2][1=]{\todo[linecolor=red,backgroundcolor=red!25,bordercolor=red,#1]{#2}}
\newcommandx{\change}[2][1=]{\todo[linecolor=blue,backgroundcolor=blue!25,bordercolor=blue,#1]{#2}}
\newcommandx{\info}[2][1=]{\todo[linecolor=OliveGreen,backgroundcolor=OliveGreen!25,bordercolor=OliveGreen,#1]{#2}}
\newcommandx{\improvement}[2][1=]{\todo[linecolor=Plum,backgroundcolor=Plum!25,bordercolor=Plum,#1]{#2}}
\newcommandx{\thiswillnotshow}[2][1=]{\todo[disable,#1]{#2}}

\DeclareMathAlphabet{\mathcalligra}{T1}{calligra}{m}{n}

\graphicspath{{Images/}}

\tikzset{  font={\fontsize{15pt}{12}\selectfont}}

\title{A massively parallel non-overlapping Schwarz preconditioner for PolyDG methods in brain electrophysiology\footnote{\textbf{Funding}: PFA is partially funded by the European Union (ERC SyG, NEMESIS, project number 101115663). Views and opinions expressed are, however, those of the authors only and do not necessarily reflect those of the European Union or the European Research Council Executive Agency. Neither the European Union nor the granting authority can be held responsible for them. CBLS has been funded by the National Recovery and Resilience Plan (NRRP), Mission 4, Component 1 – Investment 3.4 and Investment 4.1, funded by the European Union. The present research is part of the activities of the Dipartimento di Eccellenza 2023-2027 grant, funded by MUR. PFA, SP, and CBLS are members of INdAM-GNCS.}}

\author[1]{Caterina B. Leimer Saglio \orcidlink{0009-0007-7887-919X}}

\author[1]{Stefano Pagani\orcidlink{0000-0002-6662-3433}}

\author[1]{Paola F. Antonietti\orcidlink{0000-0002-2138-3878}}

\affil[1]{MOX-Dipartimento di Matematica, Politecnico di Milano, Piazza Leonardo da Vinci 32, Milan, 20133, Italy}

\begin{document}
\maketitle

\begin{abstract}
We investigate non-overlapping Schwarz preconditioners for the algebraic systems stemming from high-order discretizations of the coupled monodomain and Barreto-Cressman models, with applications to brain electrophysiology. The spatial discretization is based on a high-order Polytopal Discontinuous Galerkin (PolyDG) method, coupled with the Crank-Nicolson time discretization scheme with explicit extrapolation of the ion term. To improve solver efficiency, we consider additive Schwarz preconditioners within the PolyDG framework, which combines (massively parallel) local subdomain solvers with a coarse-grid correction. Numerical experiments demonstrate robustness with respect to the discretization parameters, as well as a significant reduction in iteration counts compared to the unpreconditioned solver. These features make the proposed approach well-suited for parallel large-scale simulations in brain electrophysiology.
\end{abstract}

\section{Introduction}\label{sec:1}
Numerical simulations of brain electrophysiology present numerous computational challenges arising from the need for accurate discretizations and efficient, scalable solvers.
Electrophysiological computational models are derived from multiscale systems of differential equations that describe the spatio-temporal evolution of the transmembrane potential in neural tissue \cite{cressman2009influence,schreiner2022simulating}. The monodomain model coupled with detailed conductance-based ionic models is widely used in this context,  since it offers a compromise between biophysical accuracy and computational tractability  \cite{schreiner2022simulating,saglio2024high}. However, its numerical approximation is characterized by significant challenges due to the sharp traveling wavefront, the strong nonlinearity introduced by the ionic currents, and the brain tissue heterogeneity and anisotropy. To cope with this complexity, realistic simulations require high-order discretizations on very complex geometries, motivating the use of high-order polytopal Discontinuous Galerkin (PolyDG) methods \cite{antonietti_hp_2013,CangianiGeorgoulisHouston_2014,cangiani_hp-version_2017}. This method naturally supports high-order approximations needed for the simulation of wave-propagation problems, such as high-frequency electrical activity 
\cite{saglio2024high}. On the other hand, it is known that Discontinuous approximations lead to larger linear systems to be solved at each time step, making the design of robust and scalable preconditioners fundamental for efficiency. 
Domain decomposition techniques, and in particular non-overlapping Schwarz methods, represent a powerful tool to address this challenge \cite{FengKarakashian_2001,Antonietti_Ayuso_2007,antonietti2011class, antonietti2020agglomeration}. Two-level Schwarz preconditioners on polytopal agglomerated meshes can significantly reduce the condition number of the resulting discrete system, ensuring, therefore, scalable algebraic solvers \cite{Antonietti_Ayuso_2007,antonietti2020agglomeration}. 

In this work, we consider the monodomain model coupled with the Barreto-Cressman ionic model \cite{cressman2009influence} discretized in space with high-order PolyDG methods and in time with the Crank-Nicolson scheme (with explicit extrapolation of the ion term). For the resulting algebraic system, we numerically investigate the performance of a two-level non-overlapping Schwarz preconditioner \cite{FengKarakashian_2001,Antonietti_Ayuso_2007,antonietti2011class, antonietti2020agglomeration}. The coarse mesh employed for the construction of the coarse correction in the preconditioner is obtained by agglomeration of the fine mesh. Numerical experiments are presented to assess the performance of the proposed preconditioner and to investigate its robustness with respect to the discretization parameters (mesh granularity and polynomial degree).

\section{The mathematical model}
\label{sec:mathematicalmodel}

\noindent We consider the monodomain equation \cite{sundnes2006computational} coupled with the Barreto-Cressman ionic model \cite{cressman2009influence}. Given an open, bounded domain $\Omega \in \mathbb{R}^d$, $(d=2,3)$ and a final time $T>0$, we introduce the transmembrane potential $u = u(\boldsymbol{x},t): \Omega \times [0,T] \rightarrow \mathbb{R}$, and the vector $\boldsymbol{y} = \boldsymbol{y}(\boldsymbol{x},t): \Omega \times [0,T] \rightarrow \mathbb{R}^n, n\ge1,$ containing the ion concentrations and gating variables of the ionic model. The coupled multiscale problem reads as follows: For any time $ t \in (0,T]$, find $u(\boldsymbol{x},t)$ and $\boldsymbol{y}=\boldsymbol{y}(\boldsymbol{x},t)$ such that:
\begin{equation}
    \label{eq:monodomain}
    \begin{dcases}
        \chi_m C_m  \frac{\partial u}{\partial t} - \nabla \cdot (\mathbf{\Sigma} \nabla u) +  \chi_m f(u,\boldsymbol{y}) = I_\mathrm{ext} & \mathrm{in} \; \Omega \times (0,T], \\
        \frac{\partial \boldsymbol{y}}{\partial t} + \boldsymbol{m}(u,\boldsymbol{y}) = \boldsymbol{0} &\mathrm{in} \; \Omega \times (0,T], \\
        \mathbf{\Sigma} \nabla u \cdot \boldsymbol{n} = 0  & \mathrm{on}\; \partial \Omega  \times (0,T],\\
        u(0) = u^0, \; \boldsymbol{y}(0) = \boldsymbol{y}^0 &\mathrm{in}\; \Omega.\\
    \end{dcases}
\end{equation}
Here, $\Omega = \Omega_G \cup \Omega_W$, being $\Omega_G$ the grey matter and $\Omega_W$ the white matter, respectively. The conductivity tensor is defined as $\boldsymbol{\Sigma}=\sigma_{l} \mathbbm{1} + (\sigma_{n} - \sigma_{l})\boldsymbol{n}\otimes\boldsymbol{n}$
and is assumed to be constant in time and piecewise constant in space. 
More precisely, for any $x \in \Omega$ let $\boldsymbol{l} = \boldsymbol{l}(\boldsymbol{x})$  be the direction of axonal fibers and let $\boldsymbol{n} = \boldsymbol{n}(\boldsymbol{x})$ be its normal vector.  Then, $\sigma_{l}=\sigma_{l}(\boldsymbol{x})$ and $\sigma_{n}=\sigma_{n}(\boldsymbol{x})$ correspond to the conductivity along the principal axonal directions and along the orthogonal direction, respectively.
In the numerical test cases, we employ a fully isotropic conductivity tensor for the grey matter tissue ($\sigma_{l}(\boldsymbol{x})=\sigma_{n}(\boldsymbol{x}) \ \forall \boldsymbol{x} \in \Omega_G$) and anisotropic conductivity for the white matter ($\sigma_{l}(\boldsymbol{x})>\sigma_{n}(\boldsymbol{x}) \ \forall \boldsymbol{x} \in \Omega_W$).
To close the system, we supplement \eqref{eq:monodomain} with homogeneous Neumann boundary conditions on $\partial \Omega$ and suitable initial conditions $u^0$ and $\boldsymbol{y}^0$. Finally, the dynamics of $\mathbf{y}$ is modeled by the Barreto-Cressman ionic model, a system of ordinary differential equations that describes a neuron's membrane potential and the interactions between intra- and extracellular ion concentrations as described in~\cite{cressman2009influence}. Specifically, the three ionic concentrations modeled are the intracellular sodium, the extracellular potassium, and the intracellular calcium, together with the respective gating variables that drive the opening and closing of ion channels. A more detailed description of the model and its parameters can be found in \cite{saglio2024high,leimerpadaptive}. 

\section{PolyDG formulation}
\label{sec:PolyDG}
We first present the PolyDG semi-discrete formulation of problem  \eqref{eq:monodomain}. Let \(\mathcal{T}_h\) be a polytopal mesh of the domain \(\Omega\), consisting of disjoint polygonal/polyhedral elements \(K\). For each element \(K\), we define its diameter as \(h_K\) and set \(h = \max_{K \in \mathcal{F}_h} h_K < 1\). 
We denote by \(\mathcal{F}_h^I\) the set of all interior faces and by \(\mathcal{F}_h^N\) the set lying on the boundary \(\partial \Omega\), where the definition of a \emph{face} \(F \in \mathcal{F}_h=\mathcal{F}_h^I\cup \mathcal{F}_h^N\) is the one given in \cite{CangianiGeorgoulisHouston_2014}.
Let $\mathbb{P}^p(\mathcal{T}_h)$ be the space of piecewise (discontinuous) polynomials of degree at most \(p \geq 1\) on each element \(K\in \mathcal{T}_h\), we set $V_h^{\text{DG}} = \mathbb{P}^p(\mathcal{T}_h).$ 
We  define the penalization parameter \(\eta : \mathcal{F}_h^I\rightarrow \mathbb{R}_+\) as follows:
\begin{equation}
    \eta = \eta(\boldsymbol{p},h,\boldsymbol{\Sigma}) = \eta_0 
        \{\boldsymbol{\Sigma}_K\}_A \dfrac{p^2}{\{h_K\}_H} \quad \text{on } F \in \mathcal{F}_h^I, 
    \label{eq:eta}
\end{equation} 
where $\boldsymbol{\Sigma}_K = \|\sqrt{\boldsymbol{\Sigma}|_K}\|^2_{L^2(K)}$, and $\{\cdot\}_A$ is the arithmetic average, and $\{\cdot\}_H$ is the harmonic average.
 We remark that $\eta_0$ should be chosen large enough to ensure stability. We introduce the following bilinear form $\mathcal{A}(\cdot,\cdot): V_h^{\text{DG}}\times V_h^{\text{DG}} \rightarrow \mathbb{R}$:
\begin{align}
    \mathcal{A}(u,v)
    &= \int_{\Omega} \boldsymbol{\Sigma}\nabla_h u \cdot \nabla_h v \, dx
       - \sum_{F \in \mathcal{F}_h^I} \int_F 
       \big( \averagel \boldsymbol{\Sigma} \nabla u \averager \cdot \jumpl v\jumpr
       + \jumpl u \jumpr \cdot \averagel \boldsymbol{\Sigma} \nabla v\averager \big)
       \, d\sigma 
       \nonumber \\
    &\quad
       + \sum_{F \in \mathcal{F}_h^I} \int_F \eta \,\jumpl u \jumpr \cdot \jumpl v \jumpr \, d\sigma
       \qquad \forall\, u,v \in V_h^{\mathrm{DG}},
    \label{eq:coer}
\end{align}
where $\nabla_h$ is the element-wise gradient and, for regular enough scalar-valued and vector-valued functions, the jump $\averagel \cdot \averager$ and average $\jumpl \cdot \jumpr$ operators are defined as in \cite{arnold2002unified}. By fixing a basis for $V_h^{\text{DG}}$, we denote by $\mathbf{A}_h$ the matrix representation of  \eqref{eq:coer}, by $\mathbf{M}_h$ the matrix representation of the $L^2$ inner product (mass matrix), and set $\mathcal{K}_h = C_m\chi_m  \mathbf{M}_h + \frac{\Delta t}{2}\, \mathbf{A}_h.$
The fully-discrete formulation is obtained by partitioning the interval \([0, T]\) into \(N_T\) sub-intervals \((t^{(k)}, t^{(k+1)}]\), each of length \(\Delta t\), such that \(t^{(k)} = k\Delta t\) for \(k = 0, \dots, N_T-1\). For time discretization, we adopt a second-order Crank-Nicolson scheme for the linear part, with the ion term discretized with a second-order explicit extrapolation, as in \cite{leimerpadaptive,saglio2024high}, leading to the following fully discrete system
\begin{equation}
\left\{
\label{eq:fully-discrete-1}
    \begin{aligned}
      & \mathcal{K}_h \mathbf{U}^{(k+1)}  =  \left(\chi_m C_m \mathbf{M}_h - \frac{\Delta t}{2}\mathbf{A}_h \right) \mathbf{U}^{(k)}  - \mathrm{\chi_m} \Delta t \mathbf{I}^{(k+1)} + \Delta t \mathbf{F}^{(k+1)} , \\
&\mathbf{M}_h\mathbf{Y}_l^{(k+1)} = \mathbf{M}_h\mathbf{Y}_l^{(k)} - \Delta t \mathbf{G}_l^{(k)}, & \forall l = 1,\cdots,n \;, \\
      &\mathbf{U}^0 = \mathbf{U}_0, \; \mathbf{Y}^0 = \mathbf{Y}_0.
    \end{aligned}
\right.
\end{equation}
where $n$ is the dimension of the Barreto-Cressman ionic model that can be efficiently solved (since $\mathbf{M}_h$ is block diagonal). On the other hand, the monodomain equation involving $\mathcal{K}_h$ has to be solved at each time step by a suitable Krylov-type solver (e.g., conjugate gradient). In the next section, we will present a non-overlapping Schwarz preconditioner to accelerate its convergence. 

\section{A massively parallel, non-overlapping Schwarz preconditioner}
In this section, following \cite{antonietti2020agglomeration}, see also \cite{Dryja_2016},  we introduce a massively parallel, two-level non-overlapping Schwarz preconditioner for the efficient solution of \eqref{eq:fully-discrete-1}.
Starting from our problem defined over fine mesh $\mathcal{T}_h$, we first agglomerate our fine-mesh elements to obtain a subdomain partition $\mathcal{T}_S$ consisting of $N$ non-overlapping subdomains such that each $\Omega_i$ is the union of some fine mesh elements $K \in \mathcal{T}_h$, and $\Omega = \bigcup_{i=1}^{N} \Omega_i$. We assume that such a decomposition is aligned with possible discontinuities of the coefficients. Notice that in the massively-parallel case that we investigate in this paper, each subdomain can consist of just one fine mesh element $K\in \mathcal{T}_h$, that is $\mathcal{T}_S=\mathcal{T}_h$. 
For any $i=1,\ldots, N$, let $V_{i}^{\mathrm{DG}}$ be the DG space defined as before, but associated only to the elements $K \in \mathcal{T}_h$ that are contained in $\Omega_i$.  We denote by $R_i^\top: V_{i}^{\mathrm{DG}} \to V_h^{\mathrm{DG}}$ the natural extension by zero operator, and by $R_i$ its adjoint (with respect to the $L^2$ inner product). \\

Next, we introduce a coarse mesh $\mathcal{T}_H$, of granularity $H>h$, still obtained based on employing agglomeration of fine-grid elements. By construction, there holds $\mathcal{T}_h \subseteq \mathcal{T}_H$.
For $1\leq q\leq p$, on the coarse mesh $\mathcal{T}_H$, we define $V_0^\text{DG}=\mathbb{P}^q(\mathcal{T}_K)$.
Let $R_0^\top : V_0^\text{DG} \to V_h^{\mathrm{DG}}$ denote the $L^2$-projection operator i.e.,
\[
\int_\Omega R_0^\top v_0 \, w_h \, dx
=
\int_\Omega v_0 \, w_h \, dx
\qquad
\forall\, w_h \in V_h^\text{DG},
\]
and let $R_0: V_h^{\mathrm{DG}} \to V_0^\text{DG}$  be its adjoint (with respect to the $L^2$ inner product). 

In the following, with a slight abuse of notation, we will denote by $R_i$, and $R_i^\top$, $i=0,\ldots, N$ both the operators defined as before, as well as their matrix representations (in the chosen basis). We define the local and coarse components of the preconditioner as
$\mathcal{K}_i= R_i \mathcal{K}_h R_i^\top$, $i=0,\ldots, N$.
The additive Schwarz preconditioner is  defined as
\begin{equation}\label{prec}
B_{\mathrm{ad}}^{-1}
=
\sum_{i=0}^{N} R_i^\top \mathcal{K}_i^{-1} R_i.
\end{equation}
The corresponding two-level additive Schwarz operator is given by $P_{\mathrm{ad}}=B_{\mathrm{ad}}^{-1}\mathcal{K}_h$.

The linear system \eqref{eq:fully-discrete-1} is then solved, at each time step, using the Preconditioned Conjugate Gradient (PCG) method, with the additive Schwarz operator $B_{\mathrm{ad}}^{-1}$ defined in Equation \eqref{prec} as preconditioner. We observe that, if the subdomain partition coincide with the fine mesh, i.e. $N=N_h$, being $N_h$ the number of elements of $\mathcal{T}_h$, then $B_{\mathrm{ad}}^{-1}$ is nothing but a Block-Jacobi preconditioner (where each block has the dimension of the elemental local approximation space) plus a (global) coarse correction. 
We remark that, in the purely stationary diffusive case, the spectral analysis of the massively-parallel, two-level additive Schwarz operator $P_{\mathrm{ad}}$ is provided in \cite{AntoniettiHoustonPennesiSuli_2020}; we also refer to  \cite{FengKarakashian_2001, Antonietti_Ayuso_2007} for the analysis in the case of nested subdomains and coarse partitions, and to \cite{antonietti2011class} for its high-order DG extension.

\section{Numerical results}
\label{sec:NumericalResults}
\begin{figure}[h]
    \centering
    \begin{subfigure}[t]{0.23\textwidth}
\makebox[0pt][l]{\hspace{-5.5em}\includegraphics[width=2\linewidth]{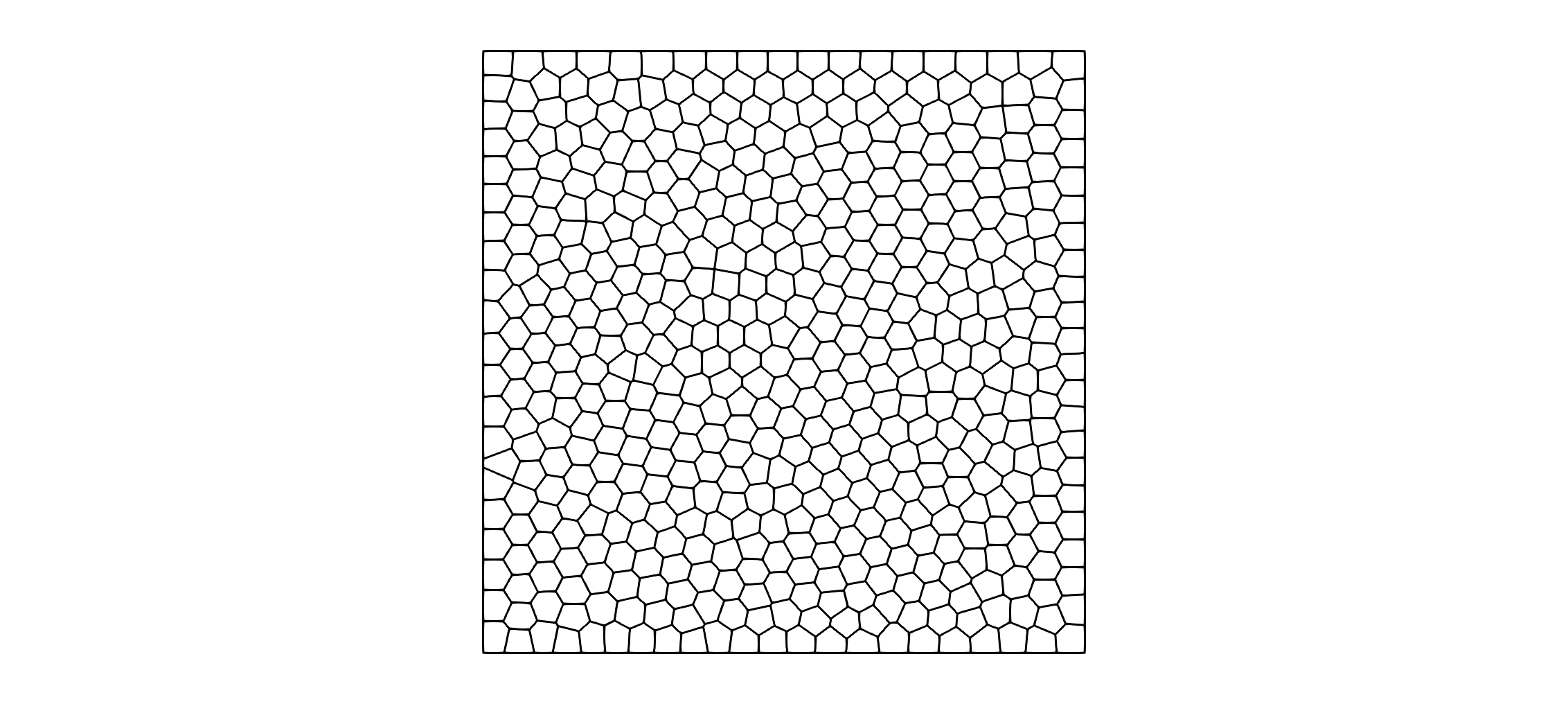}}
      \caption{Fine mesh $\mathcal{T}_h$.}
        \label{fig:1}
    \end{subfigure}
    \begin{subfigure}[t]{0.23\textwidth}
\makebox[0pt][l]{\hspace{-5.0em}\includegraphics[width=2\linewidth]{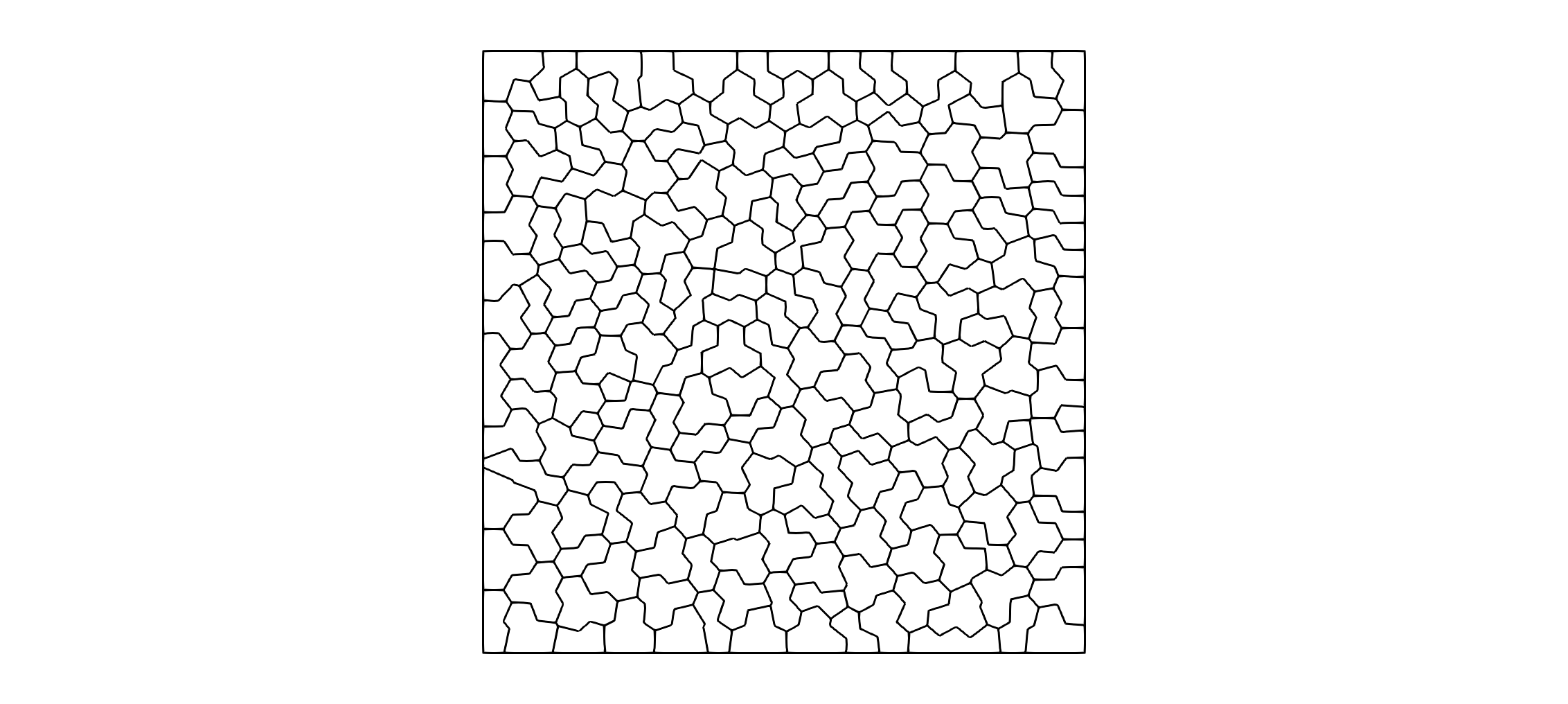}}
      \caption{$\mathcal{T}_H$ ($H = 2h$).}
        \label{fig:2}
    \end{subfigure}
        \begin{subfigure}[t]{0.23\textwidth}
\makebox[0pt][l]{\hspace{-5.0em}\includegraphics[width=2\linewidth]{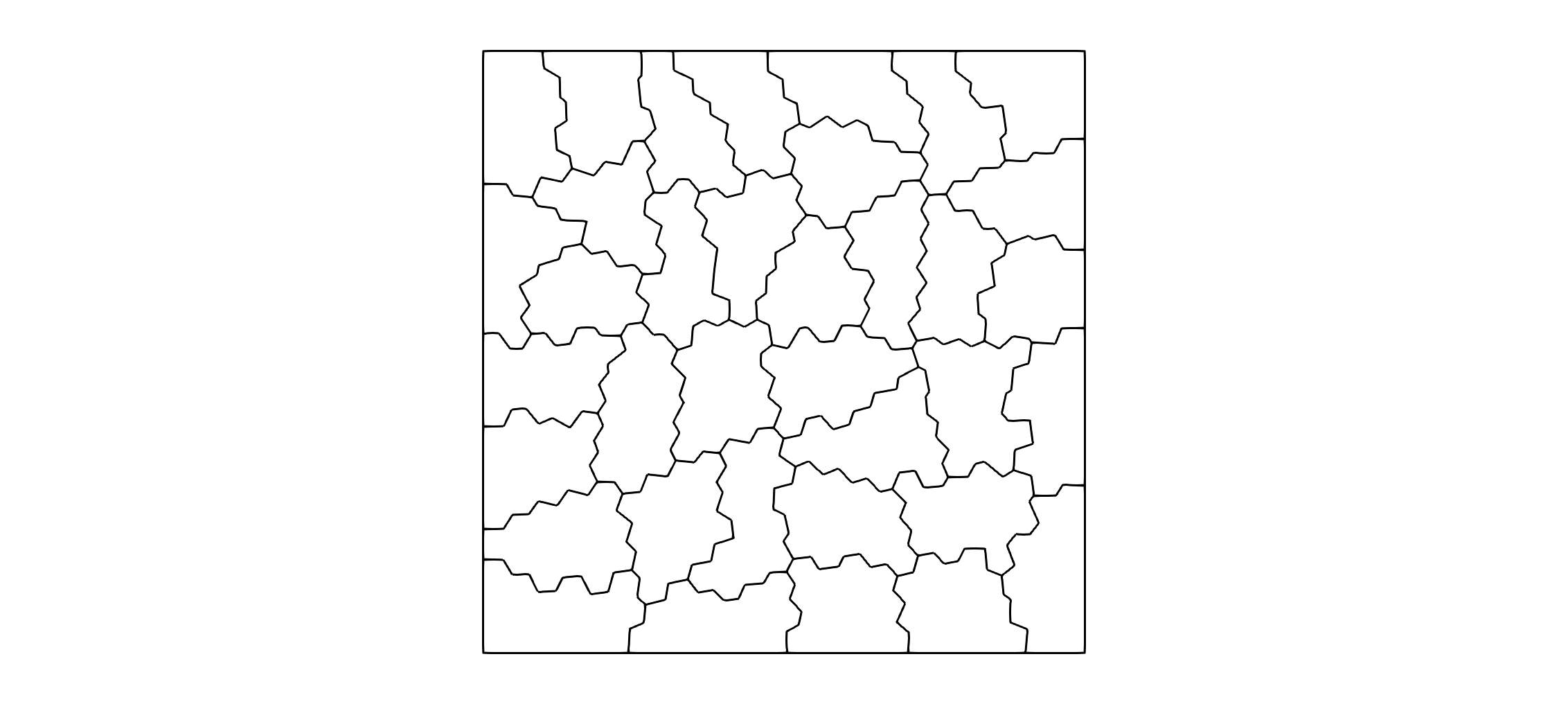}}
      \caption{$\mathcal{T}_H$ ($H = 4h$).}
        \label{fig:1}
    \end{subfigure}
        \begin{subfigure}[t]{0.23\textwidth}
\makebox[0pt][l]{\hspace{-5.0em}\includegraphics[width=2\linewidth]{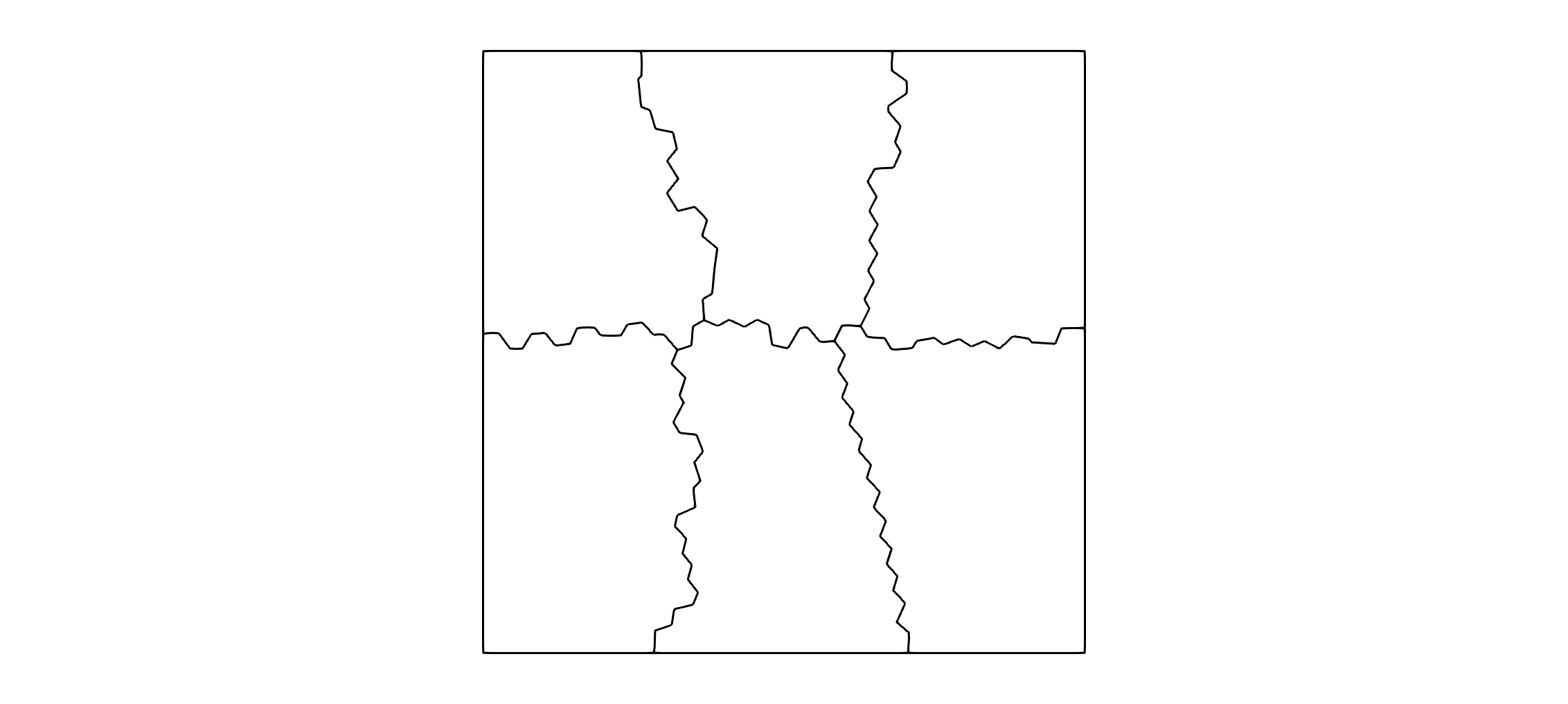}}
      \caption{$\mathcal{T}_H$ ($H = 8h$).}
        \label{fig:1}
    \end{subfigure}
        \caption{Example of a sequence of nested polygonal grids obtained by successive agglomeration of a fine mesh $\mathcal{T}_h$ consisting of 512 elements ($h \approx 0.087$).}
\label{fig:example}
\end{figure}

In this section, we present numerical experiments designed to evaluate the performance of the non-overlapping Schwarz preconditioner \eqref{prec} for the efficient solution of \eqref{eq:fully-discrete-1}. The implementation is based on the \texttt{lymph} library \cite{antonietti2025lymph}.
We investigate the robustness and scalability of the preconditioned iterative solver with respect to mesh refinement, polynomial degree, and the ratio of coarse-to-fine mesh sizes. All simulations are carried out on sequences of nested polytopal meshes, as shown in Figure~\ref{fig:example}, where we show an example of a fine mesh $\mathcal{T}_h$ consisting of 512 elements ($h \approx 0.087$) together with three examples of coarse meshes $\mathcal{T}_H$ obtained by agglomeration and with $H=2h,4h,8h$. 
We simulate the evolution of the transmembrane potential in an idealized two-dimensional square domain $\Omega$ of size $(0\,\mathrm{cm},1\,\mathrm{cm})^2$.
\begin{figure}[b]
    \centering
    \begin{subfigure}[t]{0.23\textwidth}
\makebox[0pt][l]{\hspace{-5.4em}\includegraphics[width=2\linewidth]{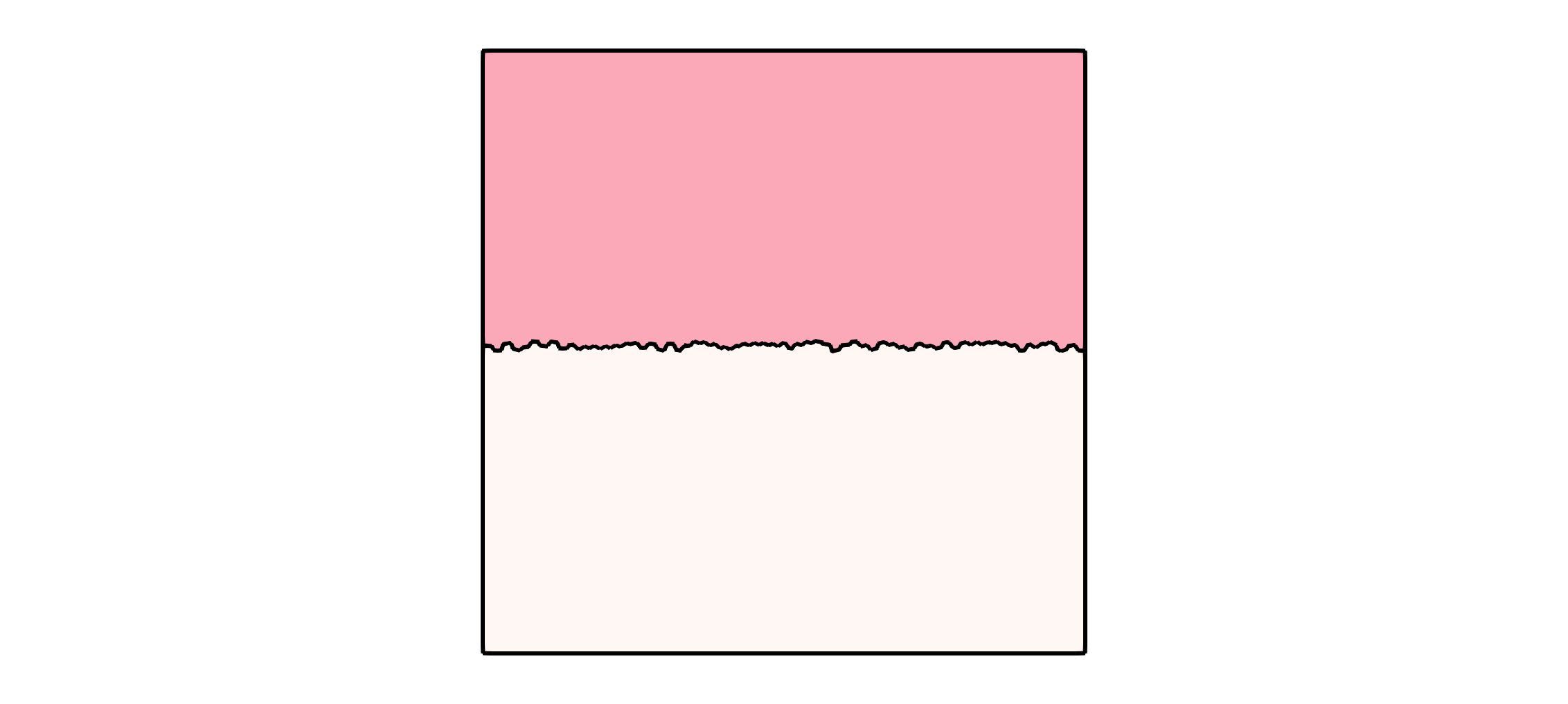}}
      \caption{$\Omega=(0,1)^2$}
        \label{fig:1}
    \end{subfigure}
    \begin{subfigure}[t]{0.23\textwidth}
\makebox[0pt][l]{\hspace{-5em}\includegraphics[width=2\linewidth]{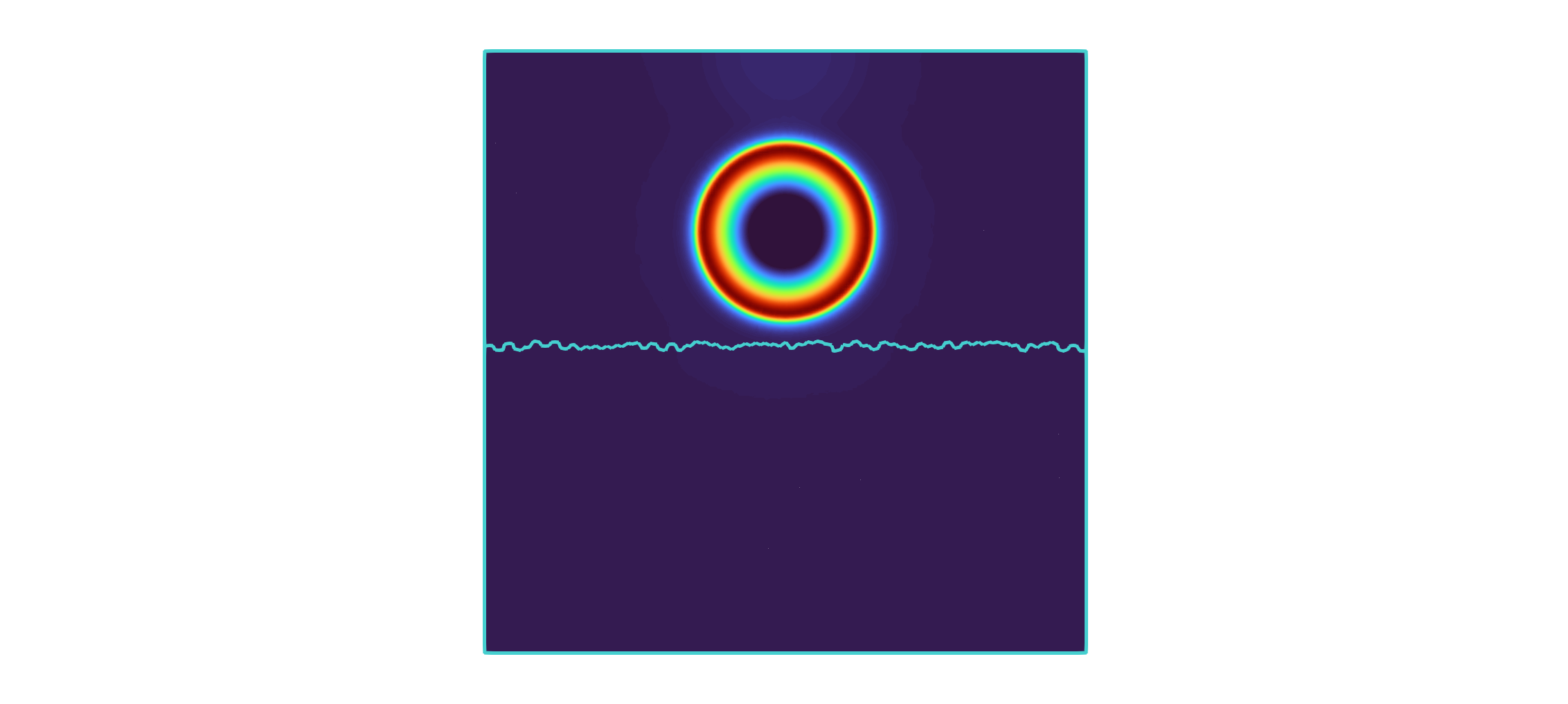}}
      \caption{ $t = 2\;\mathrm{ms}$}
        \label{fig:2}
    \end{subfigure}
        \begin{subfigure}[t]{0.23\textwidth}
\makebox[0pt][l]{\hspace{-5em}\includegraphics[width=2\linewidth]{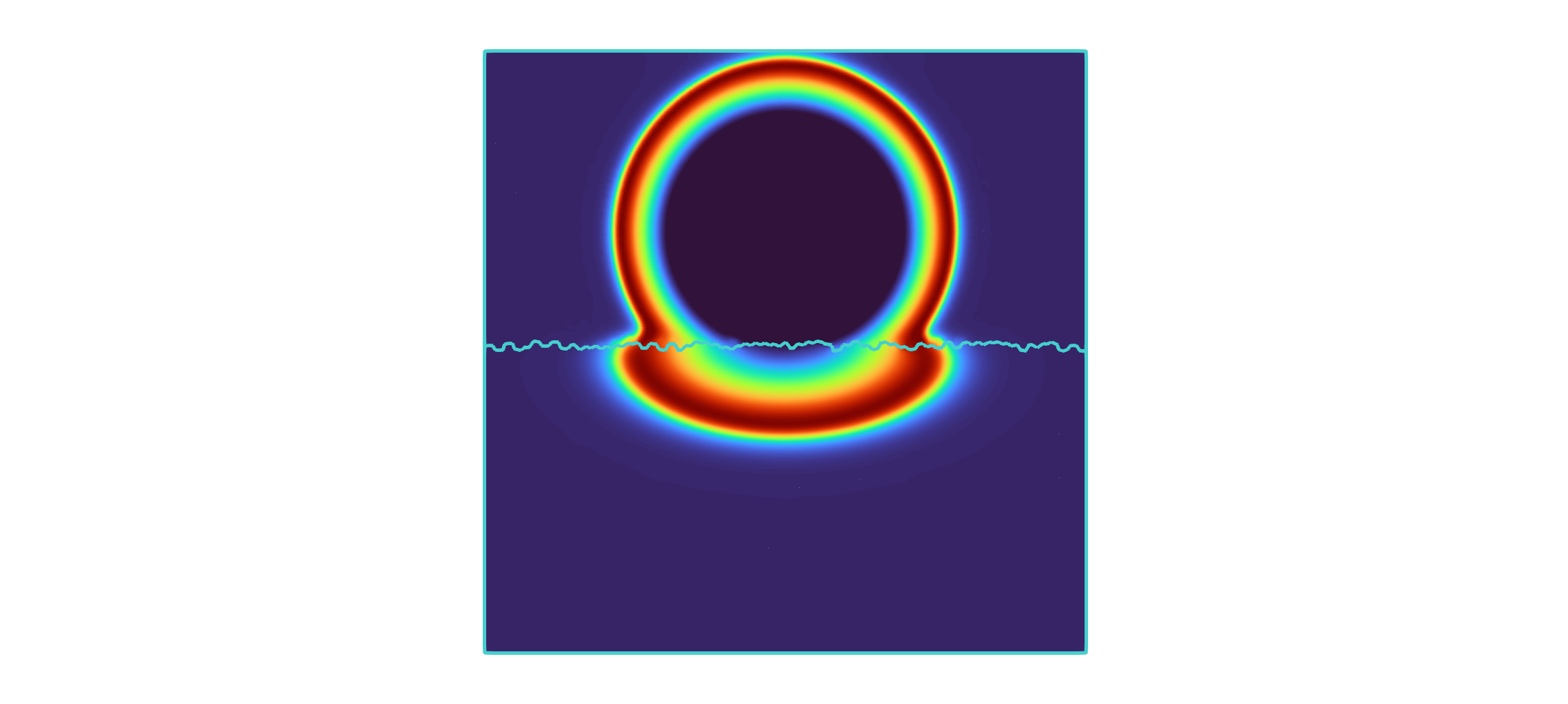}}
      \caption{$t = 4.8\;\mathrm{ms}$}
        \label{fig:1}
    \end{subfigure}
        \begin{subfigure}[t]{0.23\textwidth}
\makebox[0pt][l]{\hspace{-5em}\includegraphics[width=2\linewidth]{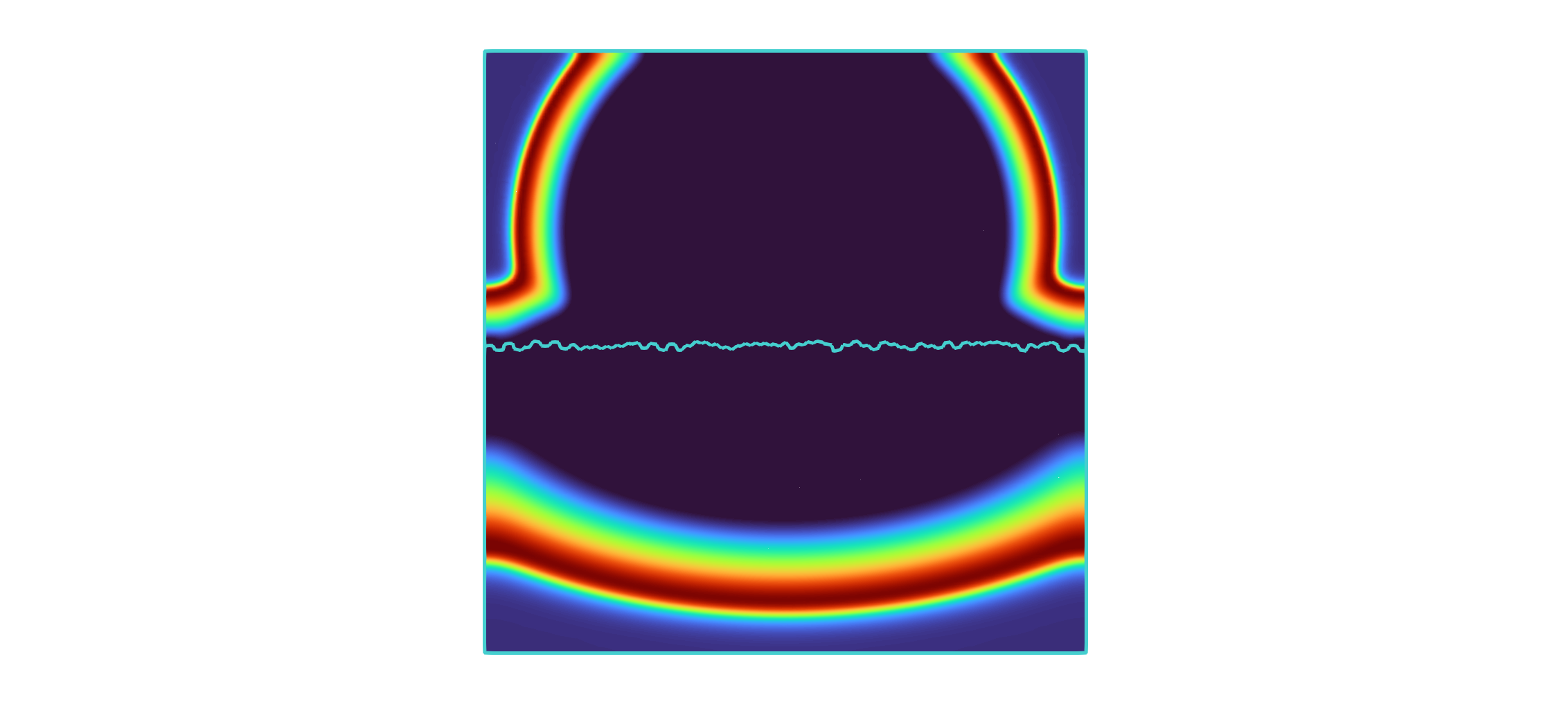}}
      \caption{$t = 8.4\;\mathrm{ms}$}
        \label{fig:1}
    \end{subfigure}
            \begin{subfigure}[t]{0.03\textwidth}
\makebox[0pt][l]{\hspace{-1.1em}\includegraphics[width=2\linewidth]{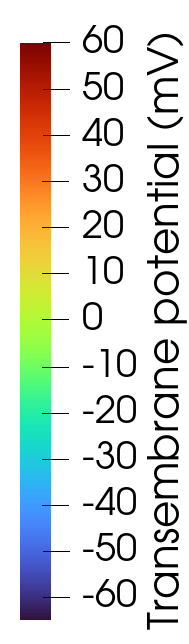}}
    \end{subfigure}
        \caption{(a) Computational domain $\Omega=(0,1)^2$: isotropic grey matter (top) and anisotropic white matter (bottom) tissue. (b)-(d) Evolution of the computed transmembrane potential for different time instants $t=2,4.8,8.4 \;\mathrm{ms}$ ($p=4$).}
\label{fig:exampleu}
\end{figure}
The domain is divided into two sub-regions, modeling grey and white matter regions, each characterized by a different conductivity value, cf. Figure~\ref{fig:exampleu}(a). In the grey matter, we set $ \sigma_l = \sigma_n = 0.63 \;\mathrm{S m ^{-1}}$ where in the white matter part we take into account anisotropy with respect to the vertical direction $\sigma_l = 0.69 \;\mathrm{S m ^{-1}}$ and $\sigma_n = 2.571 \:\mathrm{S m ^{-1}}$.  We define $\Omega_0 = \left\{ (x, y) \in \Omega \mid (x - 0.5)^2 + (y - 1)^2 < 0.016\right\}$ and we impose an initial localized potential imbalance which models a pathological brain region in $\Omega_0$ ($u^0|_{\Omega_0} = -50 \; \mathrm{mV}$). The initial value for the potential is $u^0=-67\;\mathrm{mV}$ in the remaining part of the domain.
Throughout the section we set $\Delta t = 2.5 \, \mu\mathrm{s}$ and $T=10\,\mathrm{ms}$.  We consider a sequence of polygonal fine grids with $N_h=512, 1 \, 024, 2\,048, 4\,096$ polygonal elements ($h \approx 0.087, 0.061, 0.043, 0.031$, respectively).
For each fine mesh with granularity $h$, the sequence of coarse meshes obtained by agglomeration doubles the mesh size $h$ of the fine mesh at each coarsening step, i.e., $H=2h,4h,8h,16h$.
In Figure~\ref{fig:exampleu}(b--d), we show the evolution of the approximate transmembrane potential at three time snapshots $t=2,4.8,8.4 \;\mathrm{ms}$ for the finest mesh ($N_h=4\,096$, $h\approx0.029$) and $p=4$. 
In Table~\ref{tab:1}, we report the PCG iteration counts (averaged over all time steps) needed to reduce the Euclidean norm of the relative residual below a tolerance of $10^{-9}$, when varying the size of the fine and coarse meshes. The results reported in Table~\ref{tab:1} (top) have been obtained for $p=q=1$, whereas the analogous ones obtained with $p=q=4$ are shown in Table~\ref{tab:1} (bottom). Each column/row of the table is obtained by varying $N_h$/$N_H$, respectively.
Table~\ref{tab:1} also reports the computed estimate of the condition number of $\kappa(P_{\mathrm{ad}})$ that has been obtained
by exploiting the analogies between the Lanczos technique and the PCG method.
In particular, during the PCG iteration, a tridiagonal matrix can be constructed whose extreme eigenvalues converge to those of $P_{\mathrm{ad}}$.
As a reference, Table~\ref{tab:1} also shows the corresponding quantities for the non-preconditioned CG method.
We observe that the condition number and the average number of iterations are approximately constant along the rows, indicating that the preconditioner is scalable when the ratio $H/h$ is held constant. This is consistent with what is known for non-overlapping Schwarz preconditioners for diffusion problems. Also, comparing the results with the analog ones obtained for the non-preconditioner CG method, it is clear that the preconditioner is very effective in reducing the computational burden associated with the algebraic solution.
\begin{table}[t]
\centering
\caption{Average number of iteration counts per time step and condition number estimates for the preconditioned CG solver for different fine and coarse meshes, for $p=q=1$ (top) and $p=q=4$ (bottom). Each sub-table shows in the last line the corresponding quantities for the non-preconditioned system.}
\begin{tabular}{
| >{\centering\arraybackslash}p{1.4cm}
| >{\centering\arraybackslash}p{2.8cm}
  >{\centering\arraybackslash}p{2.8cm}
   >{\centering\arraybackslash}p{2.8cm}
  >{\centering\arraybackslash}p{2.8cm}
|}
\multicolumn{5}{c}{$p = q= 1$} \\
\hline
\diagbox{$H$}{$h$}
& 0.087
& 0.061
& 0.043
& 0.031 \\
\hline
2$h$  & 55.06 (69) & 48.80 (68) &  56.6 (72) & 54.51 (69)\\
4$h$  &  167.7 (126) &  175.01 (133) & 174.18 (129) & 184.64 (134) \\
8$h$  & 792.7 (255) & 827.28 (262)  & 814.59 (259) & 1050.97 (270)\\
16$h$  & 1760.8 (335) & 1780.57 (365) & 2181.79 (448) & 2711.92 (480)\\ 
\hline
 CG   & 1.34e+4 (596) & 2.89e+4 (881) & 5.81e+4 (1256) & 1.23e+6 (1738) \\
\hline
\multicolumn{5}{c}{} \\
\multicolumn{5}{c}{$p = q= 4$} \\
\hline
\diagbox{$H$}{$h$}
& 0.087
& 0.061
& 0.043
& 0.031 \\
\hline
2$h$   & 245.84 (134) & 289.56 (136)  &  300.71 (138) & 290.38 (137)  \\
4$h$  &  863.33 (220) & 894.97 (219) & 872.96 (219) &  928.60 (229) \\
8$h$   & 3292.92 (419) & 2976.34 (390)  & 3714.64 (404)  &  4019.09 (415)\\
16$h$   & 8632.89 (695) & 8486.67 (681) & 10959.92 (736) &  13891.73 (725) \\
\hline
 CG   & 1.01e+6 (5145) & 2.2e+6 (7593) & 4.9e+6 (11296) & 9.25e+6 (15355)\\
\hline
\end{tabular}
\label{tab:1}
\end{table}

\section{Conclusion}
\label{sec:10}
We have numerically tested a two-level, massively parallel, non-overlapping Schwarz preconditioner for the linear systems arising from high-order PolyDG discretizations of a monodomain model of brain electrophysiology coupled with the Barreto-Cressman ionic dynamics.
Numerical results confirm the robustness of the proposed approach when the ratio of coarse-to-fine mesh sizes is held constant. A comparison with the non-preconditioned CG results demonstrates that the preconditioner significantly reduces the condition and the average iteration numbers. From a computational viewpoint, the preconditioner is well-suited to parallel implementations, since the local solvers involve small, independent problems that can be solved in parallel. These features make the proposed preconditioner particularly attractive for large-scale brain electrophysiology simulations.

\bibliographystyle{hieeetr}
\bibliography{sample.bib}

\end{document}